\documentclass[12pt,a4paper]{amsart}
\usepackage{amssymb,amscd}

\theoremstyle{plain}
\newtheorem{theorem}{Theorem}[subsection]
\newtheorem{lemma}[theorem]{Lemma}
\newtheorem{proposition}[theorem]{Proposition}
\newtheorem{corollary}[theorem]{Corollary}

\theoremstyle{definition}

\newtheorem{remark}[theorem]{Remark}
\newtheorem{remarks}[theorem]{Remarks}
\newtheorem{example}[theorem]{Example}

\numberwithin{equation}{subsection}

\newcommand\bA{{\mathbb A}}

\newcommand\bG{{\mathbb G}}

\newcommand\cA{{\mathcal A}}

\newcommand\aff{\operatorname{aff}}

\newcommand\diag{\operatorname{diag}}
\newcommand\id{\operatorname{id}}

\newcommand\rk{\operatorname{rk}}

\newcommand\End{\operatorname{End}}
\newcommand\GL{\operatorname{GL}}

\title{Local structure of algebraic monoids}

\author{Michel Brion}
\address{Universit\'e de Grenoble I\\
D\'epartement de Math\'ematiques\\
Institut Fourier, UMR 5582 du CNRS\\
38402 Saint-Martin d'H\`eres Cedex, France}
\email{Michel.Brion@ujf-grenoble.fr}

\begin{document}
 
\begin{abstract}
We describe the local structure of an irreducible algebraic 
monoid $M$ at an idempotent element $e$. 
When $e$ is minimal, we show that $M$ is an induced variety 
over the kernel $MeM$ (a homogeneous space) 
with fibre the two-sided stabilizer $M_e$ (a connected 
affine monoid having a zero element and a dense unit group). 
This yields the irreducibility of stabilizers and
centralizers of idempotents when $M$ is normal, 
and criteria for normality and smoothness of an arbitrary 
monoid $M$. Also, we show that $M$ is an induced variety over 
an abelian variety, with fiber a connected affine monoid
having a dense unit group. 
\end{abstract}

\maketitle

\section{Introduction}
\label{sec:introduction}

An algebraic monoid is an algebraic variety equipped with an
associative product map, which is a morphism of varieties
and admits an identity element. Algebraic monoids are closely
related to algebraic groups: the group $G$ of invertible 
elements of any irreducible algebraic monoid $M$ is a connected 
algebraic group, open in $M$. Thus, $M$ is an equivariant embedding
of its unit group $G$ with respect to the action of $G \times G$ via
left and right multiplication; this embedding has a unique closed
orbit, the kernel of the monoid.


This relationship takes a particularly precise form in the case 
of affine (or, equivalently, linear) monoids and groups. Indeed, 
by work of Vinberg and Rittatore, the affine irreducible algebraic
monoids are exactly the affine equivariant embeddings of connected
linear algebraic groups. Furthermore, any irreducible algebraic monoid
having an affine unit group is affine (see \cite{Vi95, Ri98,Ri06}).


Affine irreducible algebraic monoids have been intensively 
investigated, primarily by Putcha and Renner 
(see the books \cite{Pu88, Re05}). The idempotents play a 
fundamental role in their theory: for instance, the kernel 
contains idempotents, and these form a unique conjugacy class 
of the unit group. 


From the viewpoint of algebraic groups, the idempotents are exactly
the limit points of multiplicative one-parameter subgroups. It follows
easily that every irreducible algebraic monoid having a reductive unit
group is unit regular, that is, any element is the product of a unit
and an idempotent. Such reductive monoids are of special interest
(see the above references); their study has applications to
compactifications of reductive groups (see \cite{Ti03}) and to
degenerations of varieties with group actions (see \cite{AB04}). 


In contrast, little was known about the non-affine case until 
the recent classification of normal algebraic monoids by 
Rittatore and the author (see \cite{BR07}). Loosely speaking, 
any such monoid is induced from an abelian variety, with fibre 
a normal affine monoid. This result extends, and builds on,
Chevalley's structure theorem for connected algebraic groups (see 
\cite{Ch60,Co02}); it holds in arbitrary characteristics, 
like most of Putcha and Renner's results. More generally,
any normal equivariant embedding of a homogeneous variety
under an arbitrary algebraic group is induced from an abelian variety,
with fiber a normal equivariant embedding of a homogeneous variety
under an affine group; see \cite{Br07}, which also contains examples
showing that the normality assumption cannot be omitted.


In the present paper, we obtain a classification of 
all irreducible algebraic monoids in the spirit of \cite{BR07}: 
they are also induced from abelian varieties, but fibres are allowed
to be connected affine monoids having a dense unit group (Theorem
\ref{thm:glob}). This answers a long-standing question of Renner, see
\cite{Re84}. Also, we characterize the irreducible algebraic monoids
having a prescribed unit group $G$, as those equivariant embeddings
$X$ of $G$ such that the Albanese morphism $\alpha_X$ is affine
(Corollary \ref{cor:emb}).


Our approach differs from those of \cite{BR07,Br07}; it
relies on a local structure theorem for an irreducible algebraic
monoid $M$ at an idempotent $e$ (Theorem \ref{thm:loc}).  Loosely
speaking again, an open neighbourhood of $e$ in $M$ is an induced
variety over an open subvariety of the product $M e M$, with fibre
the two-sided stabilizer $M_e = \{ x \in M ~\vert~ x e = e x = e \}$. 
Note that $M_e$ is a closed submonoid of $M$ with the same 
identity element, and the zero element $e$; we show that 
$M_e$ is affine and connected, and its unit group $G_e$ is dense 
(Lemma \ref{lem:stab}).


When $e$ lies in the kernel, our local structure theorem takes 
a global form: the whole variety $M$ is induced over the kernel 
$M e M$, with fibre $M_e$ (Corollary \ref{cor:glob}). 
As a direct consequence, the normality or smoothness of $M$ 
is equivalent to that of $M_e$. 


This raises the question of classifying all smooth monoids having a
zero element; such a monoid is isomorphic (as a variety) to an affine
space, by Corollary \ref{cor:sm}. Another open problem arising from
our local structure theorem is the classification of all connected
algebraic monoids having a dense unit group and a zero element; to
make this problem tractable, one may assume that the unit group is
reductive.


Our results are obtained over an algebraically closed field of 
characteristic zero. They may be adapted to arbitrary 
characteristics, by considering group schemes and monoid schemes 
at appropriate places. For example, the stabilizer $M_e$ should 
be understood as a closed submonoid scheme of $M$; this subscheme 
turns out to be reduced in characteristic zero 
(Remark \ref{rem:loc2}), but this fails in positive characteristics, 
e.g., for certain non-normal affine toric varieties. 


This may be a motivation for developing a theory of monoid schemes;
note that, unlike group schemes, many monoid schemes over a field of 
characteristic zero are not reduced. For example, view the affine
plane $\bA^2$ as a monoid with product 
$(x_1,y_1) \cdot (x_2,y_2) = (x_1 x_2, x_1 y_2 + x_2 y_1)$
and unit $(1,0)$. Then the closed subscheme with ideal $(x^2, xy)$ is
a non-reduced submonoid, the affine line with a fat point at the
origin. 


This paper is organized as follows. We begin by gathering 
some basic definitions and results on algebraic varieties, 
algebraic groups and induced varieties. In Section 1, we
study various stabilizers and centralizers associated with
idempotents in affine irreducible algebraic monoids. This builds 
on work of Putcha (exposed in \cite[Chapter 6]{Pu88}), but we 
have modified some of his terminology in order to comply with
standard conventions in algebraic geometry and algebraic groups.
Section 2 is devoted to the local structure of affine irreducible
algebraic monoids, with applications to criteria for normality or 
smoothness, and to the irreducibility of stabilizers and 
centralizers in normal monoids. In the final Section 3, 
we obtain our classification theorem and derive some 
consequences, e.g., all irreducible algebraic monoids are 
quasi-projective varieties.

\bigskip

\noindent
{\bf Notation and conventions.}
We consider algebraic varieties over an algebraically closed 
field $k$ of characteristic zero; morphisms are understood to 
be $k$-morphisms. By a \emph{variety}, we mean a separated 
reduced scheme of finite type over $k$; in particular, 
varieties are not necessarily irreducible. A point will 
always mean a closed point.


An \emph{algebraic group} $G$ is a group scheme of finite 
type over $k$; then $G$ is a smooth variety, as $k$ has 
characteristic $0$. Also recall that $G$ is affine if and 
only if it is linear, i.e., isomorphic to a closed subgroup 
of some general linear group.


Given an arbitrary algebraic group $G$ and a closed subgroup 
$H \subset G$, there exists a quotient morphism $q : G \to G/H$, 
where $G/H$ is a quasi-projective variety, and $q$ is a principal 
$H$-bundle. Furthermore, $q$ is affine if and only if $H$ is 
affine. 


More generally, given a variety $Y$ where $H$ acts 
algebraically, consider the product $G \times Y$ where $H$ 
acts via $h \cdot (g,y) = (gh^{-1}, hy)$. If $Y$ is 
quasi-projective, then there exists a quotient morphism
$$
q_Y : G \times Y \to (G \times Y)/H =: G \times^H Y,
$$
where $G \times^H Y$ is a quasi-projective variety, and $q_Y$ 
is a principal $H$-bundle. The \emph{induced variety} 
$G \times^H Y$ is equipped with a $G$-action and a morphism 
$$
p_Y : G \times^H Y \to G/H
$$ 
such that $p_Y$ and $q_Y$ are $G$-equivariant; the 
(scheme-theoretic) fibre of $p_Y$ at the base point of $G/H$ is 
$H$-equivariantly isomorphic to $Y$. Moreover, $p_Y$ is affine 
if and only if $Y$ is affine (for these facts, see 
\cite[Proposition 4]{Se58a} and \cite[Proposition 7.1]{MFK94}.
In particular, $G \times^H Y$ is affine whenever $G/H$ and $Y$ 
are both affine.

\section{Stabilizers and centralizers}
\label{sec:pre}

\subsection{Stabilizers}
\label{subsec:stab}

Throughout this section, $M$ denotes an irreducible affine 
algebraic monoid, $1 \in M$ the identity element, and $G = G(M)$ 
the unit group. Then $G \times G$ acts on $M$ via 
$(x,y) \cdot z = x z y^{-1}$ 
(the \emph{two-sided action}); we also have the 
\emph{left action of $G$ on $M$} via $x \cdot y = xy$, 
and the \emph{right action} via $x \cdot y = y x^{-1}$. 

We fix an idempotent $e \in M$ and denote by
\begin{equation}\label{eqn:le}
\ell_e : M \longrightarrow M, \quad x \longmapsto e x
\end{equation}
the left multiplication by $e$. Clearly, $\ell_e$ is a retraction 
of the variety $M$ onto the closed subvariety
$e M = \{ x \in M ~\vert~ e x = x \}$.
Likewise, the right multiplication by $e$,
\begin{equation}\label{eqn:re}
r_e : M \longrightarrow M, \quad x \longmapsto x e
\end{equation}
is a retraction of $M$ onto $M e = \{ x \in M ~\vert~ x e = x \}$.

We also have a retraction of varieties
\begin{equation}\label{eqn:te}
t_e : M \longrightarrow e M e, \quad x \longmapsto e x e.
\end{equation}
Moreover, 
$e M e = \{ x \in M ~\vert~ e x e = x \} = e M \cap M e$ 
is a closed irreducible submonoid of $M$ with identity element 
$e$.

We put 
\begin{equation}\label{eqn:lstm}
M_e^{\ell} := \{ x \in M ~\vert~ x e = e \}
\end{equation}
(the \emph{left stabilizer of $e$ in $M$}) and
\begin{equation}\label{eqn:lstg}
G_e^{\ell} := G \cap M_e^{\ell}.
\end{equation}
Clearly, $M_e^{\ell}$ is a closed submonoid of $M$ with identity
element $1$ and unit group $G_e^{\ell}$. Moreover, $G_e^{\ell}$ is
dense in $M_e^{\ell}$ by \cite[Theorem 6.11]{Pu88}. 
Note that $M_e^{\ell}$ is the set-theoretic fibre of $r_e$ at $e$.
In fact, $M_e^{\ell}$ is also the scheme-theoretic fibre, 
as we shall see in Remark \ref{rem:loc1} (i).

Likewise, the \emph{right stabilizer of $e$ in $M$}, 
\begin{equation}\label{eqn:rstm}
M_e^r := \{ x \in M ~\vert~ e x = e \},
\end{equation}
is a closed submonoid of $M$ with identity element $1$, and dense 
unit group 
\begin{equation}\label{rstg}
G_e^r := G \cap M_e^r.
\end{equation}

The (two-sided) \emph{stabilizer of $e$ in $M$},
\begin{equation}\label{eqn:stm}
M_e := \{ x \in M ~\vert~ x e = e x = e \} = M_e^{\ell} \cap M_e^r,
\end{equation}
is also a closed submonoid of $M$ with identity element $1$, 
zero element $e$, and unit group 
\begin{equation}\label{eqn:stg}
G_e := G \cap M_e.
\end{equation}  
Moreover, $G_e$ is dense in $M_e$ by \cite[Theorem 6.11]{Pu88}
again.

Our notation for $M_e$ and $G_e$ differs from that of Putcha in 
\cite{Pu88}: his $M_e$ and $G_e$ are the irreducible components 
of ours that contain $1$. 

Also, note that $M_e^{\ell}$, $M_e^r$ and $M_e$ are generally 
reducible; equivalently, $G_e^{\ell}$, $G_e^r$ and $G_e$ are 
generally non-connected. This happens for many non-normal affine 
toric varieties, regarded as commutative monoids; see 
\cite[Example 6.12]{Pu88} for an explicit example. 

Yet the stabilizers in $M$ are always connected, as follows from 
the existence of a multiplicative one-parameter subgroup of $G$
with limit point $e$:

\begin{lemma}\label{lem:con}
{\rm (i)} There exists a homomorphism of algebraic groups
$\theta : \bG_m \to G_e$ 
which extends to a morphism of varieties
$\overline{\theta}: \bA^1 \to M$ such that 
$\overline{\theta}(0) = e$.
In particular, the closure of $G_e^0$ in $M$ contains $e$.

\smallskip

\noindent
{\rm (ii)} The fibres of $\ell_e$ and $r_e$ are connected.

\smallskip

\noindent
{\rm (iii)} $M_e^{\ell}$, $M_e^r$ and $M_e$ are connected.

\end{lemma}

\begin{proof}
(i) By \cite[Corollary 6.10]{Pu88}, there exists a maximal torus 
$T \subset G$ such that $e$ lies in $\overline{T}$ (the closure of 
$T$ in $M$). Since $\overline{T}$ is a (possibly non-normal) toric
variety, there exist a one-parameter subgroup 
$\theta : \bG_m \to T$ and an element $t\in T$ such that 
$\theta$ extends to a morphism 
$\overline{\theta} : \bA^1 \to M$ such that 
$\overline{\theta}(0) = t e$. 
Then $\overline{\theta}$ is a homomorphism of monoids, so that 
$\overline{\theta}(0)$ is idempotent. Since $t$ and $e$
commute, it follows that $t e = e$. Moreover, 
$\theta(x) e = \theta(x) \overline{\theta}(0) = 
\overline{\theta}(0) = e = e \theta(x)$
for all $x \in k^*$.

(ii) Clearly, $\ell_e$ is invariant under the left action of 
$G_e^{\ell}$. In particular, each fibre $\ell_e^{-1}(y)$, 
$y \in eM$, is stable under left multiplication by $G_e^0$. So, for
any $x \in \ell_e^{-1}(y)$, the closure of the orbit $G_e^0 \, x$ 
is an irreducible subvariety of $\ell_e^{-1}(y)$ containing both points 
$x$ and $e x = y$. It follows that $\ell_e^{-1}(y)$ is connected. 

(iii) The connectedness of $M_e^r$ (resp.~$M_e^{\ell}$) follows 
from (ii). To show the connectedness of $M_e$, note as above that 
the closure of any orbit $G_e^0 x$ is an irreducible subvariety 
of $M_e$ containing both points $x$ and $e x = e$.
\end{proof}

\subsection{The centralizer}
\label{subsec:cent}
Let
\begin{equation}\label{eqn:cenm}
C_M(e) := \{ x \in M ~\vert~ x e = e x \},
\end{equation}
this is the \emph{centralizer of $e$ in $M$}. Clearly, $C_M(e)$ is a
closed submonoid of $M$ with identity element $1$ and unit group
\begin{equation}\label{eqn:ceng}
C_G(e) := G \cap C_M(e).
\end{equation}
Moreover, $C_G(e)$ is connected by \cite[Theorem 6.16]{Pu88}. But
the example below (a variant of \cite[Example 6.15]{Pu88})
shows that $C_M(e)$ is generally reducible; in other words, 
$C_G(e)$ may not be dense in $C_M(e)$. 

\begin{example}\label{example}
Let $V$, $W$ be vector spaces of dimensions $m, n \geq 2$. Consider 
the multiplicative monoids $\End(V)$, $\End(W)$ and the map
$$
\varphi : \End(V) \times \End(W) \longrightarrow \End(V \otimes W),
\quad (A,B) \longmapsto A \otimes B.
$$
Then $\varphi$ is a homomorphism of monoids, and is the 
invariant-theoretical quotient by the $\bG_m$-action via 
$t \cdot (x,y) = (tx,t^{-1}y)$. Thus, the image of $\varphi$ 
is a closed normal submonoid, 
$$
M := \End(V) \otimes \End(W) \subset \End(V \otimes W).
$$
Its unit group is the quotient of $\GL(V) \times \GL(W)$ by 
$\bG_m$ embedded via $t \mapsto (t \id_V, t^{-1}\id_W)$.

Given two idempotents $e \in \End(V)$ and $f \in \End(W)$, the 
idempotent $e \otimes f \in M$ satisfies
$$
C_M(e \otimes f) = \{ x \otimes y \in M ~\vert~
xe \otimes yf = ex \otimes fy \}.
$$
It follows easily that 
$$
C_G(e \otimes f) = C_{\GL(V)}(e) \otimes C_{\GL(W)}(f),
$$ 
while
$$
C_M(e \otimes f) \supset (1-e) \otimes \End(W).
$$
Thus, $C_M(e \otimes f)$ is reducible whenever $e,f \neq 0,1$.
\end{example}

However, the centralizers in $M$ are always connected, 
as shown by the following:

\begin{lemma}\label{lem:cent}
{\rm (i)} The morphism
\begin{equation}\label{eqn:tau}
\tau_e : C_M(e) \longrightarrow e M e, \quad 
x \longmapsto x e = e x = e x e
\end{equation}
is a retraction of algebraic monoids.

\smallskip

\noindent
{\rm (ii)} The fibres of $\tau_e$ are connected. In particular,
$C_M(e)$ is connected.

\smallskip

\noindent
{\rm (iii)} We have an exact sequence of algebraic groups
\begin{equation}\label{eqn:ex}
\CD
1 @>>> G_e @>>> C_G(e) @>{\tau_e}>> G(eMe) @>>> 1.
\endCD
\end{equation}

\smallskip

\noindent
{\rm (iv)} The normalizer $N_G(G_e)$ equals $C_G(e)$. 
\end{lemma}

\begin{proof}
(i) is straightforward.

(ii) Note that $\tau_e$ is invariant under the (left or right) 
action of $G_e \subset C_G(e)$. So the assertion follows 
by arguing as in the proof of Lemma \ref{lem:con}(ii). 

(iii) Clearly, $\tau_e$ restrict to a homomorphism 
$C_G(e) \to G(eMe)$ with kernel $G_e$. This homomorphism is 
surjective by \cite[Remark 1.3(ii), Theorem 6.16]{Pu88}.

(iv) $C_G(e)$ normalizes $G_e$ by (\ref{eqn:ex}). 
Conversely, if $x \in G$ normalizes $G_e$,
then it normalizes $M_e$ (the closure of $G_e$ in $M$). Thus, $x$ 
commutes with the zero element $e$ of $M_e$.
\end{proof}

Next, we consider the action of $G$ on $M$ by conjugation. Then 
the isotropy group of $e$ is $C_G(e)$, so that the conjugacy 
class of $e$ is isomorphic to $G/C_G(e)$.

\begin{lemma}\label{lem:clo}
The $G$-conjugacy class of $e$ is closed in $M$.  In particular, 
the variety $G/C_G(e)$ is affine. 
\end{lemma}

\begin{proof}
We adapt a classical argument for the closedness of semi-simple 
conjugacy classes in affine algebraic groups. Let $B$ be a Borel 
subgroup of $G$. Since $G/B$ is complete, it suffices to 
check that the $B$-conjugacy class of $e$ is closed in $M$.
We may assume that $B$ contains a maximal torus $T$ such that 
$e \in \overline{T}$, see the proof of Lemma \ref{lem:con}. 
Then $T$ centralizes $e$, so that the $B$-conjugacy class of $e$ 
is an orbit of the unipotent radical of $B$; hence this class is
closed in the affine variety $M$, by \cite{Ro61}.
\end{proof}

\subsection{Left and right centralizers}
\label{subsec:lrc}

Let
\begin{equation}\label{eqn:celm}
C_M^{\ell}(e) := \{ x \in M ~\vert~ x e = e x e \},
\end{equation}
the \emph{left centralizer of $e$ in $M$}. 
For any $x, y \in C_M^{\ell}(e)$, we have 
\begin{equation}\label{eqn:lce}
x y e = x e y e = e x e y e = e x y e.
\end{equation}
Thus, $C_M^{\ell}(e)$ is a closed submonoid of $M$ with identity 
element $1$. The unit group of $C_M^{\ell}(e)$ equals
\begin{equation}\label{eqn:celg}
C_G^{\ell}(e) := G \cap C_M^{\ell}(e)
\end{equation}
(indeed, this is a closed submonoid of $G$, and hence a subgroup by 
\cite[3.5.1 Exercises 1 and 2]{Re05}). Moreover, $C_G^{\ell}(e)$ is 
connected by \cite[Theorem 6.16]{Pu88}. However, $C_M^{\ell}(e)$ 
is generally reducible. For instance, with the notation of Example 
\ref{example}, we have 
$C_G^{\ell}(e \otimes f) = C_G^{\ell}(e) \otimes C_G^{\ell}(f)$
while
$C_M^{\ell}(e \otimes f) \supset \End(V)(1 - e) \otimes \End(W)$.

We now extend the statement of Lemma \ref{lem:cent} to left 
centralizers: 

\begin{lemma}\label{lem:lrc}
{\rm (i)} $C_M^{\ell}(e)$ is the preimage of $e M e$ under
the morphism $r_e$ of (\ref{eqn:re}). 
Moreover, $r_e$ restricts to a retraction of algebraic monoids
\begin{equation}\label{rho}
\rho_e : C_M^{\ell}(e) \longrightarrow e M e, \quad 
x \longmapsto x e = e x e.
\end{equation}

\smallskip

\noindent
{\rm (ii)} The fibres of $\rho_e$ are connected. In particular, 
$C_M^{\ell}(e)$ is connected.

\smallskip

\noindent
{\rm (iii)} We have an exact sequence
\begin{equation}\label{eqn:ex1}
\CD
1 @>>> G_e^{\ell} @>>> C_G^{\ell}(e) @>{\rho_e}>> G(eMe) @>>> 1
\endCD
\end{equation}
and the equality
\begin{equation}\label{eqn:gcg1}
C_G^{\ell}(e) = G_e^{\ell} C_G(e).
\end{equation}

\noindent
{\rm (iv)} The normalizer $N_G(G_e^{\ell})$ equals $C_G^{\ell}(e)$. 
\end{lemma}

\begin{proof}
(i) is a direct verification.

(ii) follows from (i) together with Lemma \ref{lem:con}.

(iii) Clearly, $\rho_e$ yields a homomorphism of algebraic groups
$C_G^{\ell}(e) \to G(e M e)$ with kernel $G_e^{\ell}$. By
\ref{eqn:ex}, the restriction of this homomorphism to $C_G(e)$ is
surjective. This implies both statements.

(iv) $C_G^{\ell}(e)$ normalizes $G_e^{\ell}$ by
(\ref{eqn:ex1}). For the converse, if $x \in G$ normalizes
$G_e^{\ell}$, then $x^{-1} M_e^{\ell} x = M_e^{\ell}$ as $M_e^{\ell}$
is the closure of $G_e^{\ell}$. In particular,
$x^{-1} e x \in M_e^{\ell}$, i.e., $x^{-1} e x e = e$, and
$e x e = x e$.
\end{proof}

Next, for later use, we show that certain homogeneous spaces 
are affine:

\begin{lemma}\label{lem:aff}
{\rm (i)} The variety $G_e^{\ell}/G_e$ is isomorphic to 
$C_G^{\ell}(e)/C_G(e)$.
\smallskip

\noindent
{\rm (ii)} Both varieties $G_e^{\ell}/G_e$ and $C_G^{\ell}(e)/G_e$ 
are affine. 
\end{lemma}

\begin{proof}
(i) We have
$$ 
G_e^{\ell}/G_e = G_e^{\ell}/(G_e^{\ell} \cap C_G(e)) \cong
G_e^{\ell}C_G(e)/C_G(e) = C_G^{\ell}(e)/C_G(e),
$$
where the latter equality follows from (\ref{eqn:gcg1}).

(ii) $C_G^{\ell}(e)/C_G(e)$ is closed in $G/C_G(e)$, and hence is
affine by Lemma \ref{lem:clo}. Moreover, $C_G^{\ell}(e)/G_e$ is an
induced variety over the homogeneous space
$C_G^{\ell}(e)/C_G(e)$, with fibre $C_G(e)/G_e$. The latter is 
affine by Lemma \ref{lem:cent}; this implies the second statement.
\end{proof}

Similar assertions hold for the 
\emph{right centralizer of $e$ in $M$ resp.~$G$},
\begin{equation}\label{cer}
C_M^r(e) := \{ x \in M ~\vert~ e x  = e x e \}, \quad
C_G^r(e) := G \cap C_M^r(e).
\end{equation}
(Here again, our notation differs from that of Putcha: his 
$C_M^{\ell}(e)$ is our $C_M^r(e)$.) 
In particular, the morphism $\ell_e$ of 
(\ref{eqn:le}) yields an exact sequence of algebraic groups
\begin{equation}\label{eqn:ex2}
\CD
1 @>>> G_e^r @>>> C_G^r(e) @>{\lambda_e}>> G(eMe) @>>> 1
\endCD
\end{equation}
and the equality
\begin{equation}\label{eqn:gcg2}
C_G^r(e) = G_e^r C_G(e).
\end{equation}
 
This implies readily the following description of the
\emph{stabilizer of $e$ in $M \times M$},
$$
(M \times M)_e := \{ (x,y) \in M \times M ~\vert~ x e = e y \}.
$$
and of its stabilizer in $G \times G$, 
$$
(G \times G)_e = \{ (x,y) \in G \times G ~\vert~ x e y^{-1}= e \}.
$$
Note that $(M \times M)_e$ is a closed submonoid of the product 
monoid $M \times M$, with identity element $(1,1)$ and unit group
$(G \times G)_e$.

\begin{lemma}\label{lem:dou}
{\rm (i)} 
$$
(M \times M)_e = \{ (x,y)\in C_M^{\ell}(e) \times C_M^r(e)
~\vert~ \rho_e(x) = \lambda_e(y) \}.
$$

\noindent
{\rm (ii)} The two projections $M \times M \to M$ yield 
surjective morphisms
$(M \times M)_e \to C_M^{\ell}(e)$, $(M \times M)_e \to C_M^r(e)$
with connected fibres. In particular, $(M \times M)_e$ is connected.

\smallskip

\noindent
{\rm (iii)} The two projections $G \times G \longrightarrow G$
yield exact sequences
$$
1 \longrightarrow G_e^r \longrightarrow (G \times G)_e
\longrightarrow C_G^{\ell}(e) \longrightarrow 1, 
$$
$$
1 \longrightarrow G_e^{\ell} \longrightarrow (G \times G)_e 
\longrightarrow C_G^r(e) \longrightarrow 1.
$$
\end{lemma}

\begin{remark}\label{rem:red1}
If $G$ is reductive, then $C_G^{\ell}(e)$ and $C_G^r(e)$ are 
opposite parabolic subgroups of $G$, with common Levi subgroup 
$C_G(e)$; moreover, the unipotent radical $R_U(C_G^{\ell}(e))$ is
contained in $G_e^{\ell}$ (see \cite[Theorem 4.5]{Re05}). 

In view of (\ref{eqn:ex}), it follows that $G_e$ and $G(eMe)$ are 
reductive. Moreover, $G_e^{\ell}$ is the semi-direct product of 
$R_u(C_G^{\ell}(e))$ with $C_G(e) \cap G_e^{\ell} = G_e$. 
Likewise, $G_e^r$ is the semi-direct product of $R_u(C_G^r(e))$ with
$G_e$. 

The stabilizer $(G \times G)_e$ is described in \cite[Section 3]{AB04},
in the more general setting of stable reductive varieties. 
\end{remark}

\section{The local structure of affine irreducible monoids}
\label{sec:lsa}

\subsection{Local structure for the left action}
\label{subsec:leq}

Throughout this section, we maintain the notation and assumptions 
of Section \ref{sec:pre}. We first record the following consequence 
of a result of Putcha:

\begin{lemma}\label{lem:ope}
{\rm (i)} The product $C_G^r(e) e$ is an open affine subvariety 
of $M e$, isomorphic to $C_G^r(e)/G_e$. 

\smallskip

\noindent
{\rm (ii)} The product $C_G^r(e) G_e^{\ell}$ is an open affine 
subvariety of $G$, isomorphic to $C_G^r(e) \times^{G_e} G_e^{\ell}$ 
where $G_e$ acts on $C_G^r(e) \times G_e^{\ell}$ via 
$x \cdot (y,z) = (y x^{-1}, xz)$.
\end{lemma}

\begin{proof}
By \cite[Theorem 6.16]{Pu88}, $M e$ is contained in
$\overline{C_G^r(e)}$, the closure of $C_G^r(e)$ in $M$. Thus,
$M e \subset \overline{C_G^r(e)} e$, that is,  
$\overline{C_G^r(e)} e  = M e$. 
So $C_G^r(e) e$ is dense in $M e$. But $C_G^r(e) e$ is an orbit,
and hence is open in $M e$; the isotropy group of $e$ is
$C_G^r(e) \cap G_e^{\ell} = G_e$. Together with Lemma \ref{lem:aff},
this implies (i).

Note that
$$
C_G^r(e) e \cong C_G^r(e)/G_e = C_G^r(e)/(C_G^r(e) \cap G_e^{\ell})
\cong C_G^r(e) G_e^{\ell}/G_e^{\ell}
$$
is an open affine subvariety of $G e \cong G/G_e^{\ell}$. Since the 
morphism $r_e \vert_G : G \to G e$ is affine (as its source 
is affine), this implies (ii).
\end{proof}

Likewise, the product  $e C_G^{\ell}(e)$ is an open affine subvariety 
of $e M$, isomorphic to $C_G^{\ell}(e)/G_e$. Also, combining Lemmas
\ref{lem:lrc} and \ref{lem:ope}, we see that the product map
$G_e^r \times C_G(e) \times G_e^{\ell} \to G$
induces an isomorphism
$$
G_e^r \times^{G_e} C_G(e) \times^{G_e} G_e^{\ell} \cong
C_G^r(e) G_e^{\ell} = G_e^r C_G^{\ell}(e) 
= C_G^r(e) C_G^{\ell}(e),
$$ 
and the right-hand side is an open affine subvariety of $G$.

Next, we show that an affine neighbourhood of $e$ in $M$ is an 
induced variety relative to the left action of $C_G^r(e)$:

\begin{proposition}\label{prop:loc}
{\rm (i)} The subvariety
\begin{equation}\label{eqn:rop}
M_0^r := \{ x \in M ~\vert~ x e \in C_G^r(e) e \}
\end{equation}
is open in $M$, affine, stable under the two-sided action
of $C_G^r(e) \times C_G^{\ell}(e)$ on $M$, and contains $M_e^{\ell}$.

\smallskip

\noindent
{\rm (ii)} The product map $C_G^r(e) \times M_e^{\ell} \to M$
induces an isomorphism
\begin{equation}\label{eqn:fr}
f^r: C_G^r(e) \times^{G_e} M_e^{\ell} \longrightarrow M_0^r,
\end{equation}
equivariant under the two-sided action of the subgroup 
$C_G^r(e) \times G_e^{\ell} \subset C_G^r(e) \times C_G^{\ell}(e)$.

\smallskip

\noindent
{\rm (iii)} The scheme-theoretic intersection $M e \cap M_e^{\ell}$
consists of the (reduced) point $e$.
\end{proposition}

\begin{proof}
(i) Note that $M_0^r$ is the preimage of $C_G^r(e) e \subset Me$ 
under the morphism $r_e$ of (\ref{eqn:re}). Since that morphism 
is affine, and $C_G^r(e) e$ is open and affine (by Lemma 
\ref{lem:ope}), $M_0^r$ is open and affine as well.

Clearly, $M_0^r$ contains $M_e^{\ell}$ and is stable under 
$C_G^r(e)$. To show the stability under $C_G^{\ell}(e)$, consider
$x \in M_0^r$ and $g \in C_G^{\ell}(e)$. Then
$$
x g e = x e g e \in C_G^r(e) e C_G(e) = C_G^r(e) e,
$$
as $e g e \in e C_G(e)$ by (\ref{eqn:gcg2}). 

(ii) Since $r_e$ is equivariant under $C_G^r(e)$, the
natural map 
$$
C_G^r(e) \times^{G_e} r_e^{-1}(e) \longrightarrow 
r_e^{-1}(C_G^r(e) e) = M_0^r
$$ 
is an isomorphism, where $r_e^{-1}(e)$ denotes the scheme-theoretic
fibre. So it suffices to check the equality
\begin{equation}\label{eqn:mre}
M_e^{\ell} = r_e^{-1}(e).
\end{equation}
Clearly, $M_e^{\ell}$ is contained in $r_e^{-1}(e)$ as its maximal
closed reduced subscheme. Moreover, $M_e^{\ell}$  is stable under the
left action of $G_e$. So $C_G^r(e) \times^{G_e} M_e^{\ell}$ is a closed
subscheme of $C_G^r(e) \times^{G_e} r_e^{-1}(e)$, and both have the
same  closed points. But $C_G^r(e) \times^{G_e} r_e^{-1}(e)$ is an
open subscheme of $M$, and hence is reduced; this implies 
(\ref{eqn:mre}).

(iii) By (ii), $f^r$ restricts to an isomorphism
$$
C_G^r(e) \times^{G_e} (Me \cap M_e^{\ell}) \cong Me \cap M_0^r.
$$
Moreover, $e$ is the unique closed point of $Me \cap M_e^{\ell}$.
Since $Me \cap M_0^r$ is an irreducible variety, it follows that 
$Me \cap M_e^{\ell} = \{ e\}$ as schemes, by arguing as in the 
proof of (ii).

\end{proof}

\begin{remarks}\label{rem:loc1}
(i) As shown in the above proof, $M_e^{\ell}$ is the
scheme-theoretic fibre of $r_e$ at $e$. Also, $M_e^{\ell}$ may be
regarded as a slice at $e$ to the orbit $C_G^r(e) e$, or to its 
closure $Me$. 

\smallskip

\noindent
(ii) The right action of $C_G^{\ell}(e)$ on the open subvariety
\begin{equation}\label{eqn:lop}
M_0^{\ell} := \{ x \in M ~\vert~ e x \in e C_G^{\ell}(e) \}
\end{equation}
is described in similar terms.

\smallskip

\noindent
(iii) By the argument of Proposition \ref{prop:loc}, 
the product of $M$ induces an open immersion 
$G \times^{G_e^{\ell}} M_e^{\ell} \to M$;
this yields a local structure result for the left action of 
$G$. However, the orbit $Ge \cong G/G_e^{\ell}$ is generally 
not affine; this happens, for example, if $M = \End(V)$ and 
$e \neq 0, 1$. As a consequence, the variety 
$G \times^{G_e^{\ell}} M_e^{\ell}$ is generally not affine 
either.
\end{remarks}

\subsection{Local structure for the two-sided action}
\label{subsec:two}

We now show that an affine neighbourhood of $e$ in $M$ is an 
induced variety relative to the two-sided action of 
$C_G^r(e) \times G_e^{\ell}$.

\begin{theorem}\label{thm:loc}
{\rm (i)} The subvariety
\begin{equation}\label{eqn:op}
M_0 := \{ x \in M ~\vert~ x e \in C_G^r(e) e \text{ and }
e x \in e C_G^{\ell}(e) \} = M_0^r \cap M_0^{\ell}
\end{equation}
is open in $M$, affine, stable under the two-sided action of
$C_G^r(e) \times C_G^{\ell}(e)$ on $M$, and contains $M_e$.

\smallskip

\noindent
{\rm (ii)} The product map 
$C_G^r(e) \times M_e \times G_e^{\ell} \to M$ 
induces an isomorphism
\begin{equation}\label{eqn:f}
f : C_G^r(e) \times^{G_e} M_e \times^{G_e} G_e^{\ell} 
\longrightarrow M_0,
\end{equation}
equivariant under $C_G^r(e) \times G_e^{\ell}$.

\smallskip

\noindent
{\rm (iii)}  The scheme-theoretic intersection $M e M \cap M_e$ 
consists of the (reduced) point $e$.
\end{theorem}

\begin{proof}
(i) follows fom Proposition \ref{prop:loc} together with the fact 
that the intersection of any two affine open subvarieties is affine.

(ii) By Proposition \ref{prop:loc} again, the natural map
$$
C_G^r(e) \times^{G_e} (M_e^{\ell} \cap M_0) 
\longrightarrow M_0
$$
is an isomorphism. Thus, it suffices to show that the natural map
$$
M_e \times^{G_e} G_e^{\ell} \to  M_e^{\ell} \cap M_0 
$$
is an isomorphism. 

Let $x \in M_e^{\ell} \cap M_0$, then $x e = e$ and $e x = e g$ 
for some $g \in C_G^{\ell}(e)$. Thus, $e = e x e = e g e = g e$,
that is, $ g \in G_e^{\ell}$. Hence
$$
M_e^{\ell} \cap M_0 = \{ x \in M_e^{\ell} ~\vert~ 
e x \in e G_e^{\ell} \}.
$$
Since $G_e^{\ell}$ is open and dense in $M_e^{\ell}$, 
the product $e G_e^{\ell}  \cong G_e^{\ell}/G_e$ is open and dense 
in $e M_e^{\ell}$. Thus, the natural map
$$
(M_e^{\ell} \cap M_e^r) \times^{G_e} G_e^{\ell} \longrightarrow 
M_e^{\ell} \cap M_0
$$
is an isomorphism, where $M_e^{\ell} \cap M_e^r$ denotes the 
scheme-theoretic intersection. The latter intersection equals 
$M_e$ as a set, and hence as a scheme by the argument of 
Proposition \ref{prop:loc}. This yields the desired isomorphism.

(iii) By the argument of Proposition \ref{prop:loc} again, it suffices
to check that $M e M \cap M_e = \{ e\}$ as sets. For this, recall that 
$M$ is isomorphic to a closed submonoid of the multiplicative monoid 
$\End(V)$, where $V$ is a finite-dimensional vector space; see
\cite[Theorem 3.15]{Pu88}. Let $x \in M e M \cap M_e$. Then  
$\rk(x) \le \rk(e)$ and $x = e + y$ where $y \in \End(V)$ satisfies 
$y e = e y = 0$; thus, $\rk(x) = \rk(e) + \rk(y)$. It follows that 
$y = 0$, and $x = e$. 
\end{proof}

\begin{remarks}\label{rem:loc2} 
(i) $M_e$ may be regarded as a slice at $e$ to the orbit
$C_G^r(e) e G_e^{\ell}$, or to its closure $M e M$. Moreover,
$M_e$ (regarded as a closed subscheme of $M$) is reduced and equals 
the scheme-theoretic intersection of $M_e^{\ell}$ and $M_e^r$.

\smallskip

\noindent
(ii) One may wonder whether this local structure result extends 
to the two-sided action of the whole group $G \times G$. The answer
is positive for reductive monoids and minimal idempotents, by a
corollary of the Luna slice theorem (see \cite[Lemma 4.3]{AB04}). 

However, the answer is generally negative: if $M_0$ is a 
$G \times G$-stable neighbourhood of $e$ admitting an equivariant
morphism $f$ to the orbit $G e G \cong (G \times G)/(G \times G)_e$, 
then $M_0$ contains the open orbit $G \cong (G \times G)/\diag(G)$. 
Thus, the isotropy group, $\diag(G)$, is contained in a conjugate 
of $(G \times G)_e$ in $G \times G$. But this does not hold in
general, e.g., when $M = \End(V)$ and $e \neq 0,1$. 
\end{remarks}

The left and right actions do not play symmetric roles in the 
statement of Theorem \ref{thm:loc}. We now reformulate this result 
in a symmetric way, and apply it to the local structure of the
centralizer of $e$:

\begin{corollary}\label{cor:irc}
{\rm (i)} With the notation of Theorem \ref{thm:loc}, 
the product of $M$ induces isomorphisms
$$
C_G(e) \times^{G_e} M_e \cong C_M(e) \cap M_0
$$
and
$$
G_e^r \times^{G_e} (C_M(e) \cap M_0) \times^{G_e} G_e^{\ell} 
\cong M_0.
$$

\smallskip

\noindent
{\rm (ii)} $C_M(e) \cap M_0$ is irreducible. In particular, 
$C_M(e)$ is irreducible at $e$.
\end{corollary}

\begin{proof}
(i) Let $g \in C_G^r(e)$, $x \in M_e$ and $h \in G_e^{\ell}$ be such
that $g x h \in C_M(e)$. Then
$$
g e = g x h e = e g x h = e g e x h = e g e h
$$
so that $g e = e g e h e = e g e$. Thus, 
$g \in C_G^{\ell}(e) \cap C_G^r(e) = C_G(e)$. It follows that 
$e = e h$, that is, $h \in G_e^r \cap G_e^{\ell} = G_e$. Combined with
Theorem \ref{thm:loc}, this implies the first assertion. The second
assertion is a consequence of that theorem in view of the isomorphism 
$$
C_G^r(e) \cong G_e^r \times^{G_e} C_G(e),
$$
which follows in turn from (\ref{eqn:gcg2}).

(ii) By (i), $C_G(e) M_e$ is an open neighborhood of $e$ in $C_M(e)$. 
Moreover, $C_G(e)$ is dense in $C_G(e) M_e$, since $G_e$ is dense in 
$M_e$. Thus, $C_G(e) M_e = C_M(e) \cap M_0$ is irreducible.
\end{proof}

Similar arguments yield:

\begin{corollary}\label{cor:irr}
{\rm (i)} With the notation of Theorem \ref{thm:loc}, the product of
$M$ induces isomorphisms
$$
C_G^r(e) \times^{G_e} M_e \cong C_M^r(e) \cap M_0
$$
and
$$
(G \times G)_e \times^{G_e \times G_e} (M_e \times M_e) \cong 
(M_0 \times M_0)_e.
$$

\smallskip

\noindent
{\rm (ii)} $C_M^r(e) \cap M_0$ is irreducible. In particular, 
$C_M^r(e)$ is irreducible at $e$. 
\end{corollary}

Another geometric consequence of Proposition \ref{prop:loc} and
Theorem \ref{thm:loc} is the following normality criterion:

\begin{corollary}\label{cor:no1}
If $M$ is normal at $e$, then: 

\smallskip

\noindent
{\rm (i)} The stabilizers $M_e^{\ell}$, $M_e^r$ and $M_e$ are 
irreducible and normal. In particular, $G_e^{\ell}$, $G_e^r$ and 
$G_e$ are connected. 

\smallskip

\noindent
{\rm (ii)} The two-sided stabilizer $(G \times G)_e$ is connected 
as well.

\smallskip

\noindent
{\rm (iii)} $C_M(e)$, $C_M^{\ell}(e)$, $C_M^r(e)$ and $(M \times M)_e$
are normal at $e$.

\smallskip

Conversely, if one of the varieties $M_e^{\ell}$, $M_e^r$, $M_e$,
$C_M(e)$, $C_M^{\ell}(e)$, $C_M^r(e)$, $(M \times M)_e$ is normal at 
$e$, then $M$ is also normal at $e$.
\end{corollary}

\begin{proof}
(i) Denote by 
$$
\nu : \widetilde{M}_e^{\ell} \to M_e^{\ell}
$$ 
the normalization map of $M_e^{\ell}$. Then the left action of 
$G_e$ on $M_e^{\ell}$ lifts to an action on 
$\widetilde{M}_e^{\ell}$. This yields a finite morphism
$$ 
\varphi : C_G^r(e) \times^{G_e} \widetilde{M}_e^{\ell} 
\longrightarrow C_G^r(e) \times^{G_e} M_e^{\ell}
$$
which restricts to an isomorphism over a dense open subvariety. Since
$C_G^r(e) \times^{G_e} M_e^{\ell}$ is irreducible and normal (by the
normality of $M$ and Proposition \ref{prop:loc}), $\varphi$ is an
isomorphism. Thus, $\nu$ is an isomorphism, that is, $M_e^{\ell}$ is
normal; this variety is also connected by Lemma \ref{lem:con}, and
hence irreducible. It follows that $G_e^{\ell}$ is irreducible as
well.

The same argument shows that $M_e^r$ is normal. Likewise, the
normality of $M_e$ follows from Lemma \ref{lem:con} and Theorem
\ref{thm:loc}.

(ii) follows from the connectedness of $G_e$ in view of Lemma 
\ref{lem:dou}.

(iii) is a consequence of the normality of $M_e$ together with
Corollaries \ref{cor:irc} and \ref{cor:irr}.

The converse statement is proved similarly.
\end{proof}

\begin{remark}\label{rem:red2}
Assume that $G$ is reductive. Then the natural map
$$
R_u(C_G^r(e)) \times (C_G(e) \times^{G_e} M_e) \times 
R_u(C_G^{\ell}(e)) \longrightarrow M
$$
is an open immersion with image $M_0$. This statement follows from
Theorem \ref{thm:loc} combined with Remark \ref{rem:red1};
alternatively, this may be deduced from a local structure theorem for
actions of reductive groups, see \cite[Section 6]{Ti03} or
\cite[Lemma~2.8]{AB04}.

Also, note that each orbit of $M$ for the (left or right) $G$-action 
contains an idempotent. Hence the above statement describes the local 
structure of $M$ at an arbitrary point. 
\end{remark}

\subsection{The case of a minimal idempotent}
\label{subsec:min}

In this subsection, we assume that the idempotent $e$ is
\emph{minimal}, that is, $e$ is the unique idempotent of $e M e$;
equivalently, $e$ lies in the \emph{kernel} $\ker(M)$, the unique
closed orbit of $G \times G$ in $M$. Hence
\begin{equation}\label{eqn:ker}
\ker(M) = G e G = M e M.
\end{equation}
Moreover, $e M e$ is an algebraic group with identity element 
$e$; the $G$-conjugates of $e$ are exactly the minimal idempotents 
of $M$ (for these results, see \cite[Chapter 6]{Pu88} and 
\cite[Section 1]{Hu05}). Combined with (\ref{eqn:ex}), it follows that 
\begin{equation}\label{eqn:sma}
e M e = e G e = e C_G(e) = C_G(e) e.
\end{equation}
Furthermore, 
\begin{equation}\label{eqn:grl}
G = C_G^r(e) C_G^{\ell}(e)
\end{equation}
by \cite[Theorem 6.30 and Corollary 6.34]{Pu88}. In view of 
Lemma \ref{lem:lrc} (iii), this implies in turn:
\begin{equation}\label{eqn:gd}
G = C_G^r(e) G_e^{\ell} = G_e^r C_G^{\ell}(e).
\end{equation}
We now show that the open subvarieties that occur in Proposition 
\ref{prop:loc} and Theorem \ref{thm:loc} are all equal to $M$:

\begin{lemma}\label{lem:min}
{\rm (i)} $M e = C_G^r(e) e$; equivalently, $M = C_G^r(e) M_e^{\ell}$.
Likewise, $e M = e C_G^{\ell}(e)$ and $M = M_e^r C_G^{\ell}(e)$.

\smallskip

\noindent
{\rm (ii)} $M_e^r e = G_e^r e$; equivalently, 
$M_e^r = G_e^r M_e$. Likewise, 
$e M_e^{\ell} = e G_e^{\ell}$ and $M_e^{\ell} = M_e G_e^{\ell}$.

\smallskip

\noindent
{\rm (iii)} $M = C_G^r(e) M_e G_e^{\ell} = C_G^{\ell}(e) M_e G_e^r$
and $\ker(M) = C_G^r(e) e G_e^{\ell} = C_G^{\ell}(e) e G_e^r$.

\smallskip

\noindent
{\rm (iv)} $M_0^r = M_0^{\ell} = M_0 = M$.
\end{lemma}

\begin{proof}
(i) By (\ref{eqn:ker}), (\ref{eqn:sma}) and (\ref{eqn:gd}),
$M e = M e M e = G e G e = G e C_G(e) = G e = 
C_G^r(e) G_e^{\ell} e =  C_G^r(e) e$.

(ii) Let $x \in M_e^r$, then $x e \in C_G^r(e) e = C_G(e) G_e^r e$.
Write accordingly $x e = g h e$, then 
$$
e = e x = e x e = e g h e = g e h e = g e.
$$
Thus, $g \in G_e$ and $x \in G_e^r e$. 

(iii) follows from (i) and (ii) together with (\ref{eqn:ker});
likewise, (iv) follows from (i) and (iii).
\end{proof}

Together with Theorem \ref{thm:loc} and Corollaries \ref{cor:irc} and
\ref{cor:irr}, this lemma implies the following global structure 
result:

\begin{corollary}\label{cor:glob}
For any minimal idempotent $e$, the product of $M$ induces
isomorphisms
$$
C_G^r(e) \times^{G_e} M_e \times^{G_e} G_e^{\ell} \cong M,
\quad  M_e \times^{G_e} G_e^{\ell} \cong M_e^{\ell},
$$
$$
C_G^r(e) \times^{G_e} M_e \cong C_M^r(e), \quad
C_G(e) \times^{G_e} M_e \cong C_M(e),
$$
$$
\text{and} \quad
(G \times G)_e \times^{G_e \times G_e} (M_e \times M_e) 
\cong (M \times M)_e.
$$
\end{corollary}

Also, $M$ is normal if and only if it is normal at some minimal 
idempotent, since the normal locus is stable under the 
two-sided $G \times G$-action. Together with Corollaries 
\ref{cor:no1} and \ref{cor:glob}, this implies in turn:

\begin{corollary}\label{cor:nor}
Let $e$ be a minimal idempotent of an irreducible algebraic monoid
$M$. Then the following assertions are equivalent:

\smallskip

\noindent
{\rm (i)} $M$ is normal.

\smallskip

\noindent
{\rm (ii)} All the varieties $M_e$, $M_e^{\ell}$, $M_e^r$, 
$C_M(e)$, $C_M^{\ell}(e)$, $C_M^r(e)$ and $(M \times M)_e$ 
are irreducible and normal.

\smallskip

\noindent
{\rm (iii)} At least one of these varieties is normal at $e$.
\end{corollary}

\section{The structure of irreducible monoids}

\subsection{Local structure}
\label{subsec:ls}

In this subsection, we extend most results of the previous
sections to an arbitrary (possibly non-affine) irreducible 
algebraic monoid $M$ with unit group $G$. 

As in \cite{BR07} which treats the case where $M$ is normal, our main
tool is a theorem of Chevalley: there exists a unique exact sequence
of connected algebraic groups 
\begin{equation}\label{eqn:chev}
1 \longrightarrow G_{\aff} \longrightarrow G \longrightarrow 
\cA(G) \to 0,
\end{equation}
where $G_{\aff}$ is affine and $\cA(G)$ is an abelian variety
(see \cite{Ch60}, and \cite{Co02} for a modern proof).
It follows that $G_{\aff}$ is the maximal closed connected affine
subgroup of $G$, while the quotient morphism 
$$
\alpha_G : G \longrightarrow \cA(G)
$$ 
is the \emph{Albanese morphism of the variety $G$}
(the universal morphism to an abelian variety, see \cite{Se58b}).

Denote by $M_{\aff}$ the closure of $G_{\aff}$ in $M$. Clearly,
$M_{\aff}$ is a submonoid of $M$ with identity element $1$ and unit
group $G_{\aff}$. In fact, $M_{\aff}$ is affine by 
\cite[Theorem 2]{Ri06}; as a consequence, $M_{\aff}$ is the maximal
closed irreducible affine submonoid of $M$. Moreover, the natural map
\begin{equation}\label{eqn:mor1}
\pi : G \times^{G_{\aff}} M_{\aff} \longrightarrow M, \quad 
(g,x) G_{\aff} \longmapsto g x
\end{equation}
is birational (since $\pi$ restricts to an isomorphism
$G \times^{G_{\aff}} G_{\aff} \to G$) and proper (since 
$G/G_{\aff} \cong \cA(G)$ is complete). It follows that 
$\pi$ is surjective, that is,
\begin{equation}\label{eqn:mgm}
M = G M_{\aff}.
\end{equation}

Let $C$ denote the centre of $G$; then $G = C G_{\aff}$ (see e.g.
\cite[Lemme 2]{Se58a}). As a consequence,
\begin{equation}\label{eqn:gzg}
G = C^0 G_{\aff} \quad \text{and} \quad M = C^0 M_{\aff},
\end{equation}
where $C^0$ denotes the neutral component of $C$. In particular, 
\begin{equation}\label{eqn:alb}
C^0/(C^0 \cap G_{\aff}) \cong G/G_{\aff} \cong \cA(G),
\end{equation}
and the natural map 
\begin{equation}\label{eqn:mor2}
\pi : C^0 \times^{C^0 \cap G_{\aff}} M_{\aff} \longrightarrow M
\end{equation}
is proper and birational. This yields the following generalization 
of \cite[Corollary 2.4]{BR07}: 

\begin{lemma}\label{lem:ide}
Any idempotent of $M$ is contained in $M_{\aff}$.
\end{lemma}

\begin{proof}
Given $x \in M$, the (set-theoretical) fibre of $\pi$ at $x$ may be
identified with the subvariety
$$
\{ z(C^0 \cap G_{\aff}) ~\vert~ z \in C^0, ~ z^{-1} x \in M_{\aff} \}
\subset C^0/(C^0 \cap G_{\aff}) \cong \cA(G).
$$
If $x$ is idempotent, then the above subvariety is a closed subsemigroup
of $\cA(G)$, and hence is a group by 
\cite[3.5.1 Exercises 1 and 2]{Re05}. It follows that $x \in M_{\aff}$.
\end{proof}

We now choose an idempotent $e \in M_{\aff}$ and define the stabilizers 
$M_e^{\ell}, M_e^r, M_e \subset M$ and $G_e^{\ell}, G_e^r, G_e \subset G$
as in Section \ref{subsec:stab}. Then again, $M_e^{\ell}$ is a
submonoid of $M$ with identity element $1$ and unit group
$G_e^{\ell}$, and likewise for $M_e^r$, $M_e$. 

\begin{lemma}\label{lem:afis}
The stabilizers $G_e^{\ell}$, $G_e^r$ and $G_e$ are affine.
\end{lemma}

\begin{proof}
Recall that $G_e^{\ell}$ is the isotropy group of the point 
$e \in M$ for the left $G$-action. Since this action is faithful, 
$G_e^{\ell}$ is affine by \cite[Lemma p.~54]{Ma63}. So
$G_e^r$ and $G_e = G_e^{\ell} \cap G_e^r$ are affine as well.
\end{proof}

\begin{remarks}
(i) The two-sided stabilizer $(G \times G)_e$ is not necessarily 
affine, as it contains $C^0$ embedded diagonally in $G \times G$.

\smallskip

\noindent
(ii) In general, the stabilizers are not contained in $G_{\aff}$, 
as shown by \cite[Example 2.7]{BR07}. Specifically, let $A$ be 
a non-trivial abelian variety, $F \subset A$ a non-trivial finite 
subgroup, and $M$ the commutative monoid obtained from the 
product monoid $A \times \bA^1$ by identifying the points $(x,0)$ 
and $(x + f,0)$, for all $x \in A$ and $f \in F$. Then  
$G = A \times \bG_m$, $G_{\aff} = \bG_m$, and the image of 
$(0,0)$ in $M$ is an idempotent with stabilizer
$$
G_e = G_e^{\ell} = G_e^r = F \times \bG_m,
$$
which strictly contains $G_{\aff}$.
\end{remarks}

Similary, we may define the centralizers
$C_M^{\ell}(e), C_M^r(e), C_M(e) \subset M$ and 
$C_G^{\ell}(e), C_G^r(e), C_G(e) \subset G$ as in Section
\ref{subsec:cent}. Then 
$$
C_M^{\ell}(e) = C^0 C_{M_{\aff}}^{\ell}(e), 
\quad
C_M^r(e) = C^0 C_{M_{\aff}}^r(e),
\quad 
C_M(e) = C^0 C_{M_{\aff}}(e)
$$
and likewise
$$
C_G^{\ell}(e) = C^0 C_{G_{\aff}}^{\ell}(e), 
\quad
C_G^r(e) = C^0 C_{G_{\aff}}^r(e),
\quad 
C_G(e) = C^0 C_{G_{\aff}}(e).
$$
In particular, these closed subvarieties are all connected, and
$$
C_G^r(e) = G_e^r C_G(e), \quad C_G^{\ell}(e) = G_e^{\ell} C_G(e).
$$
Moreover, $M e = C^0 G_{\aff} e$ contains 
$C^0 C^r_{G_{\aff}}(e) e = C_G^r(e) e$
as a dense open subvariety, by Lemma \ref{lem:ope}. 

Thus, \emph{all the statements of Proposition \ref{prop:loc} and
Theorem \ref{thm:loc} hold in this setting, except for the affineness
of $M_0^r$, $M_0^{\ell}$ and $M_0$}; the proofs are exactly the same. 
\emph{Corollaries \ref{cor:irc}, \ref{cor:irr} and \ref{cor:no1} 
hold as well}, since their proofs do not use any assumption of 
affineness. Combined with the following result, this reduces the local 
structure of irreducible algebraic monoids to that of connected affine
monoids having a dense unit group.

\begin{lemma}\label{lem:stab}
The stabilizers $M_e^{\ell}$, $M_e^r$ and $M_e$ are affine and connected.
Their unit groups $G_e^{\ell}$, $G_e^r$, $G_e$ are dense. 
\end{lemma}

\begin{proof}
With the notation of Proposition \ref{prop:loc}, the preimage
$$
(f^r)^{-1}(G \cap M_0^r) = C_G^r(e) \times^{G_e} (G \cap M_e^{\ell})
=  C_G^r(e) \times^{G_e} G_e^{\ell}
$$
is dense in $C_G^r(e) \times^{G_e} M_e^{\ell}$, as $G$ is dense in $M$.
It follows that $G_e^{\ell}$ is dense in $M_e^{\ell}$. 
Since $G_e^{\ell}$ is affine, this implies the affineness of $M_e^{\ell}$
in view of \cite[Theorem 3]{Ri06}. 

The connectedness of $M_e^{\ell}$ is obtained by arguing as in the 
proof of Lemma \ref{lem:con}.

Likewise, the desired properties of $M_e$ follow from the statement
of Theorem \ref{thm:loc}.
\end{proof}

These considerations yield the following smoothness criterion:

\begin{corollary}\label{cor:sm}
Let $e$ be an idempotent of an irreducible algebraic monoid $M$. Then
$e$ is a smooth point of $M$ if and only if the variety $M_e$ 
is an affine space.
\end{corollary}

\begin{proof}
By Theorem \ref{thm:loc}, $M$ is smooth at $e$ if and only if $M_e$ is 
smooth at $e$. So the assertion follows from the existence of an 
attractive $\bG_m$-action on $M_e$ with fixed point $e$. 

Specifically, let $\theta : \bG_m \to G_e$ be as in Lemma
\ref{lem:con}. Then the $\bG_m$-action on $G_e$ via 
$t \cdot x = \theta(t) x$ extends to an action of the multiplicative 
monoid $\bA^1$ on $M_e$, such that $0 \cdot x = e x = e$ for all 
$x \in M_e$. This yields a positive grading of the algebra of regular 
functions on the affine variety $M_e$. Now the graded version of 
Nakayama's lemma implies our assertion.
\end{proof}

\begin{remark}
The above smoothness criterion raises the question of classifying  
algebraic monoid structures on a given affine $n$-space,
having the origin as their zero element. When the unit group is 
reductive, such structures correspond bijectively to decompositions 
of $n$ into a sum of squares of positive integers, 
as the corresponding monoids are just products of matrix monoids.

Indeed, if $M$ is a smooth monoid with reductive unit group 
$G$ and zero element $0$, then the variety $M$ is equivariantly 
isomorphic to the $G \times G$-module $T_0 M$, as follows from 
the Luna slice theorem. This yields a $G \times G$-equivariant
isomorphism $\varphi : M \to \prod_{i=1}^m \End(V_i)$, 
where $V_1,\ldots,V_m$ are simple $G$-modules; as a consequence,
$\varphi$ is an isomorphism of monoids.
Thus, $G$ is identified to an open subgroup of the product
$\prod_{i=1}^m \GL(V_i)$, and hence to the whole product. 
(This is also proved in \cite[Section 11]{Ti03}, via a 
representation-theoretic argument.)
\end{remark}

In the case that $e$ is a minimal idempotent of $M_{\aff}$, 
the subvariety $G_{\aff} e G_{\aff}$ is the unique closed orbit of 
$G_{\aff} \times G_{\aff}$ in $M_{\aff}$. As $\pi$ is proper, 
it follows that $G e G = C^0 G_{\aff} e G_{\aff}$ 
is the unique closed $G \times G$-orbit in $M$. In other words, 
$G e G$ is the kernel of $M$. Then \emph{all the statements of 
Section 2.3 hold}, with exactly the same proofs. 

Also, note that the minimal idempotents of $M$ are exactly those 
of $M_{\aff}$ (by Lemma \ref{lem:ide}); they form a unique 
conjugacy class of $G_{\aff}$ or, equivalently, of $G$ by 
(\ref{eqn:gzg}). 

Since the smooth locus of $M$ is stable under the two-sided 
action of $G \times G$, we see that \emph{$M$ is smooth if and 
only if $M_e$ is an affine space for some minimal idempotent $e$.}

\subsection{Global structure} 
\label{subsec:gls}

By the main result of \cite{BR07}, the map $\pi$ of (\ref{eqn:mor1}) 
is an isomorphism whenever $M$ is \emph{normal}, and then
$M_{\aff}$ is normal as well. In other words, any normal monoid is 
an induced variety over an abelian variety, with fibre a normal affine 
monoid. 

This statement does not extend to arbitrary irreducible monoids, in
view of \cite[Example 2.7]{BR07}. Yet we show that any such monoid
is an induced variety over an abelian variety, with fibre a connected 
affine monoid having a dense unit group:

\begin{theorem}\label{thm:glob}
Let $M$ be an irreducible algebraic monoid, and $G$ its unit group. 
Then there exists a closed submonoid $N \subset M$ satisfying 
the following properties:  

\smallskip

\noindent
{\rm (i)} $N$ is affine, connected, and contains $1$.

\smallskip

\noindent
{\rm (ii)} The unit group $H := G(N)$ is dense in $N$, and contains 
$G_{\aff}$ as a subgroup of finite index. In particular, $M_{\aff}$
is the irreducible component of $N$ containing $1$.

\smallskip

\noindent
{\rm (iii)} The canonical map 
\begin{equation}\label{eqn:iso}
\varphi : G \times^H N \longrightarrow M, \quad (g,n)H \longmapsto gn
\end{equation}
is an isomorphism of varieties.

\smallskip

Moreover, the projection $p : G \times^H N \to G/H$ is identified 
with the Albanese morphism of the variety $M$. In particular,
$H$ and $N$ are uniquely determined by $M$. 
\end{theorem}

\begin{proof} 
We begin with the proof of the final assertion: we assume that
$M = G \times^H N$ where $H$ and $N$ satisfy (i)--(iii), and
show that $p$ equals the Albanese morphism 
$\alpha_M : M \to \cA(M)$. The latter morphism is uniquely determined 
up to a translation in $\cA(M)$; we normalize it by imposing
that $\alpha_M(1) = 0$ (the origin of the abelian variety 
$\cA(M)$).

Consider a morphism (of varieties)
$$
\alpha : G \times^H N \to A
$$
where $A$ is an abelian variety. The restriction of $\alpha$ 
to the neutral component $H^0 \subset N$ is a morphism from 
a connected affine algebraic group to an abelian variety, and 
hence is constant (as follows e.g. from 
\cite[Corollary 3.9]{Mi86}. Thus, $\alpha$ is 
constant on every irreducible component of $N$. Since $N$ is 
connected, $\alpha$ maps $N$ to a point; likewise, it maps each 
fibre of $p$ (that is, each translate $gN$ in $G \times^H N$) 
to a point. Together with Zariski's Main Theorem, this implies 
that $\alpha$ factors as $p$ followed by a morphism 
$G/H \to A$. This proves the desired equality.

In particular, $N$ is identified with the fibre of $\alpha_M$ at 
$0$. We now show that this fibre satisfies the properties 
(i)--(iii). 

By rigidity, the restriction $\alpha_M \vert_G$  is a 
homomorphism of algebraic groups (see e.g. 
\cite[Corollary 3.6]{Mi86}). Thus, $\alpha_M$ is a homomorphism
of algebraic monoids. In particular, $N$ is a closed submonoid 
of $M$ containing $1$.

Moreover, $\alpha_M \vert_G$ factors through a unique 
homomorphism $\cA(G) \to \cA(M)$, which is surjective as $G$ 
is dense in $M$. Since $\cA(G) = G/G_{\aff}$, we may identify 
$\cA(M)$ with the homogeneous space $G/H$, where $H$ is a 
closed subgroup of $G$ containing $G_{\aff}$. This identifies 
$M$ with $G \times^H N$, equivariantly for the right $G$-action
on $M$.  

Since $M$ is connected, it follows that $H$ acts transitively 
on the connected components of $N$. Let $N'\subset N$ 
be the connected component containing $1$, and $H'\subset H$
its stabilizer. Then the canonical map $H \times^{H'} N' \to N$
is an isomorphism, as $H/H'$ is identified with the set of 
connected components of $N$. Thus, the analogous map
$G \times^{H'} N' \to M$ is an isomorphism as well. Moreover, 
since $H'$ has finite index in $H$, and $G/H$ is complete, 
it follows that $G/H'$ is complete as well. Thus, the composite 
map $M \cong G \times^{H'} N' \to G/H'$ factors through a 
$G$-equivariant morphism $G/H \to G/H'$. This implies 
that $H = H'$ and $N = N'$, i.e., $N$ is connected. 

Likewise, since $G$ is dense in $M$, it follows that $H$ is 
dense in $N$. To complete the proof, it suffices to show that
the quotient $H/G_{\aff}$ is finite. Indeed, this implies that 
$H$ is affine and, in turn, that $N$ is affine in view of
\cite[Theorem 3]{Ri06}. 

The finiteness of $H/G_{\aff}$ is equivalent to the assertion
that the canonical homomorphism 
$$
G/G_{\aff} = \cA(G) \to \cA(M) = G/H
$$ 
has a finite kernel, and hence to the existence of a 
$G$-equivariant morphism 
$$
\psi: M \to \cA(G)/F,
$$ 
where $F \subset \cA(G)$ is a finite subgroup.

To construct such a morphism $\psi$, choose a minimal idempotent
$e \in M$ and recall that $M e = G e \cong G/G_e^{\ell}$ 
(see Lemma \ref{lem:min}). This yields a $G$-equivariant 
morphism $\gamma : M \longrightarrow G/G_e^{\ell}$. Now let
$$
\psi :  M \longrightarrow \cA(G/G_e^{\ell})
$$
be the composition of $\gamma$ with the Albanese morphism of
$G/G_e^{\ell}$. Then $\psi$ is $G$-equivariant. Moreover, 
$\cA(G/G_e^{\ell})$ is the quotient of $\cA(G)$ by the image 
of the subgroup $G_e^{\ell}$, and the latter image is a finite group 
(as $G_e^{\ell}$ is affine by Lemma \ref{lem:afis}). 
\end{proof}

\begin{remark}\label{rem:ind}
We may define a natural structure of algebraic monoid on 
$G \times^H N$ so that the map $\varphi$ of (\ref{eqn:iso}) 
is an isomorphism of algebraic monoids. Indeed, the canonical map
$$
C^0 \times^{C^0 \cap H} N \to G \times^H N
$$
is an isomorphism, as $G = C^0 H \cong C^0 \times^{C^0 \cap H} H$.
Moreover, $C^0 \times^{C^0 \cap H} N$ is the quotient of the product
monoid $C^0 \times N$ by the central subgroup $C^0 \cap H$,
embedded via $x \mapsto (x,x^{-1})$. 

Alternatively, one may observe that the $H$-action on $N$ by
conjugation extends uniquely to a $G$-action, where $C^0$ acts
trivially (since $C^0 \cap H$, a central subgroup of $H$, acts 
trivially on $N$ by conjugation). Thus, one may form the
semi-direct product of monoids $G \times N$: its product is given by
$$
(x,a) \cdot (y,b) = (xy, a^{y^{-1}}b)
$$
where $a^z$ denotes the conjugate of $a \in N$ by $z \in G$
(see \cite[Example 3.7]{Re05}). Then one checks that this product
induces a unique product on $G \times^H N$ such that the quotient 
map $G \times N \to G\times^H N$ is a homomorphism of monoids.

The above construction is an analogue for algebraic monoids of the
induction of varieties with group actions. 
\end{remark}

\subsection{Some applications}
\label{subsec:sa}

We begin by stating two direct consequences of Theorem \ref{thm:glob}, 
first obtained in \cite{BR07} for normal monoids:

\begin{corollary}\label{cor:qp}
Any irreducible algebraic monoid is quasi-projective.
\end{corollary}

\begin{corollary}\label{cor:equiv}
The category of irreducible algebraic monoids is 
equi\-valent to the category having as objects the triples $(G,H,N)$, 
where $G$ is a connected algebraic group, $H \subset G$ is a closed 
subgroup containing $G_{\aff}$ as a subgroup of finite index, 
and $N$ is a connected affine algebraic monoid with unit group $H$, 
dense in $N$.

The morphisms from such a triple $(G,H,N)$ to a triple $(G',H',N')$
are the pairs $(\varphi, \psi)$, where $\varphi : G \to G'$ is a
homomorphism of algebraic groups such that $\varphi(H) \subset H'$,
and $\psi : N \to N'$ is a homomorphism of algebraic monoids such that 
$\varphi \vert_H =  \psi \vert_H$.
\end{corollary}

Another consequence is a characterization of monoids among (possibly
non-normal) equivariant embeddings of algebraic groups:

\begin{corollary}\label{cor:emb}
Let $G$ be a connected algebraic group and let $X$ be a 
$G \times G$-equivariant embedding of $G$. Then $X$ admits a (unique) 
structure of algebraic monoid if and only if $\alpha_X$ is affine. 

\end{corollary}

\begin{proof} 
If $X$ is an irreducible algebraic monoid, then its Albanese morphism
is affine by Theorem \ref{thm:glob}. For the converse, arguing as in 
the proof of that theorem, one shows that $\cA(X) \cong G/H$, where 
$H \subset G$ is a closed subgroup containing $G_{\aff}$; moreover, 
$\alpha_X$ is $G$-equivariant. Thus, $X \cong G \times^H Y$, 
where $Y$ is an equivariant embedding of the (possibly non-connected) 
algebraic group $H$. Moreover, $Y$ is affine by assumption, and 
hence is an algebraic monoid. In particular, its unit group $H$ is 
affine, so that $H/G_{\aff}$ is finite. As in Remark \ref{rem:ind}, 
the induced variety $X$ is then an algebraic monoid.
\end{proof}

Next, we show how to recover the main result of \cite{BR07}
(Theorem 4.1 and its proof):

\begin{corollary}\label{cor:fin}
For any irreducible algebraic monoid $M$, the morphism 
$\pi : G \times^{G_{\aff}} M_{\aff} \to M$ of (\ref{eqn:mor1}) 
is finite.

In particular, $M$ is normal if and only if the associated triple 
satisfies: $H = G_{\aff}$ and $N = M_{\aff}$ is normal.
\end{corollary}
  
\begin{proof}
Since $\pi$ is proper and $G$-equivariant, and $M = GN$, 
the finiteness of $\pi$ is equivalent to the finiteness
of its restriction to $\pi^{-1}(N)$. But 
$$
\pi^{-1}(N) = H \times^{G_{\aff}} M_{\aff}
$$ 
by Theorem \ref{thm:glob}. Furthermore, 
$\pi \vert_{\pi^{-1}(N)}$ factors as the closed embedding
$$
H \times^{G_{\aff}} M_{\aff} \to H \times^{G_{\aff}} N
$$
(corresponding to the inclusion of $M_{\aff}$ into $N$), followed by
the isomorphism
$$
H \times^{G_{\aff}} N \cong (H/G_{\aff}) \times N
$$
(since $N$ is $H$-stable), followed in turn by the projection
$$
(H/G_{\aff}) \times N \to N,
$$
a finite morphism.

Since $\pi$ is birational and finite, it is an isomorphism 
whenever $M$ is normal, by Zariski's Main Theorem; it then follows 
that $M_{\aff}$ is normal as well. Moreover, $H = G_{\aff}$ and 
$N = M_{\aff}$ by the uniqueness statement in 
Theorem \ref{thm:glob}. The converse is obvious.
\end{proof}

Finally, one may show as in \cite[Theorem 5.3]{BR07} that 
\emph{any irreducible algebraic monoid $M$ has a faithful 
representation by endomorphisms of a homogeneous vector bundle 
over an abelian variety (the Albanese variety of $M$.)}


\begin{thebibliography}{100} 




\bibitem[AB04]{AB04}
V.~Alexeev and M.~Brion,
\emph{Stable reductive varieties I: Affine varieties},
Invent. math. {\bf 157} (2004), 227--274.


\bibitem[Br07]{Br07} 
M.~Brion,
\emph{Some basic results on actions of non-affine algebraic 
groups}, arXiv: math.AG/0702518.


\bibitem[BR07]{BR07} M.~Brion and A.~Rittatore,
\emph{The structure of normal algebraic monoids},
Semigroup Forum {\bf 74} (2007), 410--422.


\bibitem[Ch60]{Ch60} C.~Chevalley,
\emph{Une d\'emonstration d'un th\'eor\`eme sur les groupes 
alg\'ebriques}, 
J. Math. Pures Appl. (9) {\bf 39} (1960), 307--317.


\bibitem[Co02]{Co02} B.~Conrad,
\emph{A modern proof of Chevalley's theorem on algebraic groups},
J. Ramanujam Math. Soc. {\bf 17} (2002), 1--18.


\bibitem[Hu05]{Hu05} W.~Huang,
\emph{The kernel of a linear algebraic semigroup},
Forum Math. \textbf{17} (2005), 851--869.


\bibitem[Ma63]{Ma63} H.~Matsumura,
\emph{On algebraic groups of birational transformations},
Atti Accad.~Naz.~Lincei Rend.~Cl.~Sci.~Fis.~Mat.~Natur. (8)
{\bf 34} (1963), 151--155.


\bibitem[Mi86]{Mi86} J.~S.~Milne,
\emph{Abelian Varieties}, in:
Arithmetic Geometry (G.~Cornell and J.~H.~Silverman, eds.), 103--150,
Springer-Verlag, New York, 1986. 


\bibitem[MFK94]{MFK94}
D.~Mumford, J.~Fogarty and F.~Kirwan,
\emph{Geometric Invariant Theory}, 3rd enlarged edition, 
Ergeb. Math. {\bf 36}, Springer-Verlag, 1994.


\bibitem[Pu88]{Pu88} M.~S.~Putcha,
\emph{Linear Algebraic Monoids},
London Math. Soc. Lecture Note Series {\bf 133},
Cambridge University Press, Cambridge, 1988.

\bibitem[Re84]{Re84} L.~E.~Renner,
\emph{Quasi-affine algebraic monoids},
Semigroup Forum {\bf 30} (1984), 167--176. 

\bibitem[Re05]{Re05} L.~E.~Renner,
\emph{Linear Algebraic Monoids},
Invariant Theory and Algebraic Transformation Groups, V,
Encyclop\ae dia Math. Sci. {\bf 134},
Springer-Verlag, Berlin, 2005.  


\bibitem[Ri98]{Ri98} A.~Rittatore,
\emph{Algebraic monoids and group embeddings},
Transform. Groups {\bf 3} (1998), 375--396.


\bibitem[Ri06]{Ri06} A.~Rittatore,
\emph{Algebraic monoids with affine unit group are affine}, 
Transform. Groups {\bf 12} (2007), 601--605.


\bibitem[Ro61]{Ro61} M.~Rosenlicht, 
\emph{On quotient varieties and the affine embedding of certain
homogeneous spaces},
Trans. Amer. Math. Soc. {\bf 101} (1961), 211--223. 


\bibitem[Se58a]{Se58a} J.-P.~Serre, 
\emph{Espaces fibr\'es alg\'ebriques,}
S\'eminaire C.~Chevalley (1958), Expos\'e No.~1, 
Documents Math\'ematiques {\bf 1}, Soc. Math. France, Paris, 2001.


\bibitem[Se58b]{Se58b} J.-P.~Serre, 
\emph{Morphismes universels et vari\'et\'e d'Albanese},
S\'eminaire Chevalley (1958--1959), Expos\'e No. 10, 
Documents Math\'ematiques {\bf 1}, Soc. Math. France, Paris, 2001.


\bibitem[Ti03]{Ti03} D.~A.~Timashev, 
\emph{Equivariant compactifications of reductive groups},
Russ. Acad. Sci. Sb. Math. {\bf 194} (2003), No. 4, 589--616.






\bibitem[Vi95]{Vi95}
E.B.~Vinberg,
\emph{On reductive algebraic semigroups}, in:
Lie groups and Lie algebras: E.B.~Dynkin's seminar, 
Amer. Math. Soc. Transl. Ser. 2 {\bf 169} (1995), 
145--182.


\end{thebibliography}
\end{document}